\documentclass[twoside,amstex]{article}
\usepackage{amsmath}
\usepackage{upref}
\usepackage{amssymb}
\usepackage{amscd}
\def\bc{\begin{center}}
\def\ec{\end{center}}
\def\no{\noindent}

\usepackage[all]{xy}
\usepackage{graphicx,tikz}
\usepackage{amsmath,amsxtra,amssymb,latexsym, amscd,amsthm}
\usepackage{indentfirst}
\usepackage[mathscr]{eucal}
\usepackage[pagebackref=true]{hyperref}

\DeclareMathOperator{\depth}{depth}
\DeclareMathOperator{\pd}{pd}

\begin{document}
\baselineskip12pt \thispagestyle{empty} \vspace*{-2mm}
\no\hspace*{97mm}
%\begin{picture}(27.7,16)(0,0)\thicklines\setlength{\unitlength}{1mm}
%\put(0,0.2){\line(1,0){28.7}}\put(0,16){\line(1,0){28.7}}
%\put(0,12.5){\bf\slshape Algebra}\put(0,8.5){\bf\slshape Colloquium}
%\put(0,5){\footnotesize \copyright\hspace{.1mm} 2009 AMSS CAS}
%\put(0,1.5){\footnotesize \hspace*{3.2mm} \& SUZHOU UNIV}
%\end{picture}

%\vspace*{-17.6mm}\no{\small{\sl Algebra Colloquium} {\bf 1?}\,:\,?
%(200?) ??--??} \vspace*{46mm}

%\pagestyle{myheadings} \markboth
%{\hfill {\small\sl \.  Hac\i o\~glu}}%×÷ÕßËõд
%{{\small\sl weak F-regular and Strongly F-regular}\hfill}%üÌâ
%\vspace*{-1.5 true cm}

\begin{center} {\large\bf Depth of the 3-path Ideal of  Square of a Path} \end{center}
\vskip1mm

\begin{center} { { Liuqing Yang,\,\,  Lizhong Chu }} \end{center}

\vskip1mm

\no\small{{\bf Abstract.}\ \
 In this note,  we compute depth of the 3-path ideal of square of a path and  show that the 3-path ideal $I_3(P_n^2)$ of square of a path
 graph is Cohen-Macaulay if and only if $n=3$ or $4$. Also, we consider the limit behavior of depth of powers of the 3-path ideal of square of a path.

\begin{figure}[b]
\vspace{-4mm}
\rule[-2.5truemm]{5cm}{0.1truemm}\\[2mm]
{\small
  Address correspondence to Lizhong Chu, Department of
Mathematics, Soochow University,
Suzhou 215006, China; E-mail: Chulizhong@suda.edu.cn.}%»ù½ð×ÊÖú
\end{figure}

\vskip2mm\vskip2mm\vskip2mm\vskip2mm

\no{\bf 1.\ \ Introduction} \vskip2mm
\noindent
Harary and Ross (\cite{HR}) defined the squares of trees (the $2^{th}$ powers of trees) and then their definition has been extended to squares of graphs. For a finite simple graph $G$, its square, denoted by $G^2$, is the graph with the same vertex set as $G$ and two vertices are adjacent in $G^2$ if they are adjacent in $G$ or their distance in $G$ is 2.  This notion is  generalized to the $k^{th}$ power of a graph $G$ later (See \cite{Go}\cite{Lin}).
Properties of  squares of graphs  have been intensively studied in combinatorics (\cite{AA,AMH,FH,FK,LT,LT1,L,R,R1,SW}). Classes of graphs which are closed to taking squares and more general taking powers, have been determined. Strongly chordal graphs (\cite{L,R1}), interval graphs (\cite{R}), proper interval graphs (\cite{R}) are known to have this property. Moreover, many researchers paid attention to recognition of squares of graphs and they developed algorithms for their recognition. The complexity of the problem of recognition of squares of graphs or of square roots of graphs has been determined for several classes of graphs.

Given a finite simple graph $G$ on the vertex set $V(G)=\{1,\ldots,n\}$ and with the set of edges $E(G)$,  the edge ideal of $G$ in the polynomial ring $S=K[x_1,\ldots, x_n]$ over a field $K$, denoted by $I(G)$,  is defined as $$I(G)=\langle x_ix_j: \{i,j\}\in E(G)\rangle.$$ During the past decades, researchers studied algebraic properties and  homological invariants of the graph $G$ (See \cite{EV, FVT2, Ha, MT,  SM1, Vi1} for example). In particular, Z. Iqbal, M. Ishaq and M. Aamir \cite{IA} considered the Stanley Conjecture for the edge ideal of the square of a path or a cycle.   Olteanu \cite{O} studided the behaviour of algebraic and homological invariants of the edge ideal of the square of a tree.  He computed the depth and dimension of the edge ideal of the square of a path.
%From his result, it follows that the edge ideal of the square of a path on $n$ vertices is Cohen Macaulay if and only if $n=2,3$.

The notion of path ideals has been introduced in \cite{CN} as a generalization of the edge ideals. For a graph $G$, $t$-path ideal of $G$ is defined as $$I_t(G)=(x_1x_2\cdots x_t: {\text{where}\, } P: x_1, x_2, \cdots , x_t \, \text{is\, a\, } t{\text{-path\,  in}}\, G) .$$  Researchers also have been trying
to explore the algebraic invariants of $t$-path ideal (or $3$-path ideal for simple) via the combinatorial invariants of graphs in the last few decades ( See \cite{AL1, AL2, BC1, BC2, B01, B02,  DC, Ku} for example).

In this note, we aim at computing the depth of the 3-path ideal of  square of a path and determining its Cohen-Macaulay property.

\vskip2mm\vskip2mm

\setcounter{page}{1} \vskip2mm \vskip4mm
%\normalsize
%\def\Ker{{\rm Ker}}
%\def\Im{{\rm Im}}
%\def\Ind{{\rm Ind}}
%\def\rank{{\rm rank}}

\no{\bf 2.\ \ Main Results}
\vskip2mm
 \no
For two positive integers $n, m$ $(m< n)$, let ${{S}_{m}}=k[{{x}_{1}},...,{{x}_{m}}]$, ${{S}_{(m)}}=k[{{x}_{m+1}},...,{{x}_{n}}]$ and $S_n=S$.
We denote a path on $n$ vertices $\{1,2, \cdots, n\}$ by $P_n$.
For $n\geq 2$, the $2^{th}$ power of a path or square of a path, denoted by $P_n^2$, is a graph such that for all $1\leq i<j\leq n$, $\{i,j\}\in E(P^{2}_{n})$ if and only if $0<j-i\leq 2$. Then 3-path ideal $I_3(P_n^2)$ of $P_n^2$ is:
$$I_3(P_n^2)=(x_1x_2x_3, x_2x_3x_4, \cdots, x_{n-4}x_{n-3}x_{n-2}, x_{n-3}x_{n-2}x_{n-1}, x_{n-2}x_{n-1}x_{n}, $$
$$\qquad \qquad x_1x_3x_4, x_2x_4x_5, \cdots,  x_{n-5}x_{n-3}x_{n-2}, x_{n-4}x_{n-2}x_{n-1},x_{n-3}x_{n-1}x_{n},$$
$$\qquad \qquad x_1x_2x_4, x_2x_3x_5, \cdots, x_{n-5}x_{n-4}x_{n-2}, x_{n-4}x_{n-3}x_{n-1}, x_{n-3}x_{n-2}x_{n},$$
$$\qquad \qquad x_1x_3x_5, x_2x_4x_6, \cdots,  x_{n-6}x_{n-4}x_{n-2}, x_{n-5}x_{n-3}x_{n-1}, x_{n-4}x_{n-2}x_{n}).$$

 Firstly we give some lemmas. \vskip2mm

 \no{\bf Lemma 2.1} (\cite[Theorem 4.3]{GC} )\,  Let $I$ be a monomial ideal of $I$ and let $f$ be an arbitrary monomial in $S$.  Then

 (i) \, $ \depth {{S}}/I \in \{ \depth S/(I:f), \depth S/(I, f)\}.$

(ii) \,$ \depth {{S}}/I = \depth S/(I:f)$ if $\depth S/(I, f)\geq \depth S/(I:f)$.

\vskip2mm

 \no{\bf Lemma 2.2} (\cite[Theorem 3.1.34, Proposition 3.1.33]{Vi}) \, Let $I$ be an ideal of $R=K[x_1,x_2,\cdots, x_s ]$ and Let $J$ be an ideal of $T=K[y_1,y_2,\cdots, y_t ]$. Let $Q=R\otimes_K T$. Then $$\depth Q/(I+J)=\depth R/I+\depth T/J.$$

\vskip2mm

 \no{\bf Lemma 2.3}\,   For a positive integer $n$,  the following inequality holds:  $$\left\lceil \frac{n-1}{7} \right\rceil+ \left\lfloor \frac{n-3}{7} \right\rfloor  +1  \geq \left\lceil \frac{n}{7} \right\rceil+ \left\lfloor \frac{n-2}{7} \right\rfloor $$

\vskip2mm

\no{\bf Proof.}  Assume that $n=7k+r, \, 0\leq r\leq 6.$ We prove the inequality in the two cases.

In case that $2\leq r \leq 6$:

$$\left\lceil \frac{n-1}{7} \right\rceil+ \left\lfloor \frac{n-3}{7} \right\rfloor  +1 =\left\lceil \frac{7k+r-1}{7} \right\rceil+ \left\lfloor \frac{n-3}{7} \right\rfloor  +1=k+1+\left\lfloor \frac{n+4}{7} \right\rfloor ,$$
While,
$$ \left\lceil \frac{n}{7} \right\rceil+ \left\lfloor \frac{n-2}{7} \right\rfloor = k+1+\left\lfloor \frac{n-2}{7} \right\rfloor .$$

In case that $0\leq r \leq 1$:
$$\left\lceil \frac{n-1}{7} \right\rceil+ \left\lfloor \frac{n-3}{7} \right\rfloor  +1 =\left\lceil \frac{n+6}{7} \right\rceil+ \left\lfloor \frac{7k+r-3}{7} \right\rfloor  =\left\lceil \frac{n+6}{7} \right\rceil+ k-1,$$
While,
$$ \left\lceil \frac{n}{7} \right\rceil+ \left\lfloor \frac{n-2}{7} \right\rfloor = \left\lceil \frac{n}{7} \right\rceil+ k-1.$$
Then the result holds.

\vskip2mm

In the following, we compute depth of 3-path ideals of  square path graphs. \vskip2mm

\no{\bf Theorem 2.4} \, Let $I_3(P_n^2)$ be  the 3-path ideal of square of a path $P_n^2$. Then
$$\depth({{S}}/I_3(P_{n}^{2}))=\left\lceil \frac{n}{7} \right\rceil+ \left\lfloor \frac{n-2}{7} \right\rfloor  +1 .$$
\vskip2mm

\no{\bf  Proof.} Our proof mainly depends on Lemma 2.1. For $n=3,4,5,6,7$, we can verify the result's correctness by the following method with $I_3(P_{i}^2)=0, x_i=0 $ if $i<1$.  Note that
$$I_3(P_n^2):x_{n-2}=I_3(P_{n-5}^2)+ J,$$
where, $$ J=(x_{n-4}x_{n-3}, x_{n-3}x_{n-1}, x_{n-1}x_{n}, x_{n-5}x_{n-3}, x_{n-4}x_{n-1}, x_{n-5}x_{n-4}, x_{n-3}x_{n}, x_{n-6}x_{n-4}, x_{n-4}x_{n}).$$
Firstly, we compute $ \depth S/((I_3(P_n^2):x_{n-2}):x_{n-4})$.
$$((I_3(P_n^2):x_{n-2}):x_{n-4})=I_3(P_{n-5}^2)+( x_{n-6},  x_{n-5}, x_{n-3},  x_{n-1}, x_{n})$$
$$\qquad \qquad \qquad \qquad \qquad \quad =I_3(P_{n-7}^2)+( x_{n-6},  x_{n-5}, x_{n-3},  x_{n-1}, x_{n}),$$
By Lemma 2.2, we have that
$$ \depth S/((I_3(P_n^2):x_{n-2}):x_{n-4})=\depth S/(I_3(P_{n-7}^2)+( x_{n-6},  x_{n-5}, x_{n-3},  x_{n-1}, x_{n}))  \qquad   \qquad  $$ $$  \qquad   \qquad  \qquad   \qquad  \qquad   \qquad  \qquad  \, =\depth S_{n-7}/I_3(P_{n-7}^2)+ \depth S_{(n-7)}/( x_{n-6},  x_{n-5}, x_{n-3},  x_{n-1}, x_{n})   $$ $$\quad  \, \, =\depth S_{n-7}/I_3(P_{n-7}^2)+ 2.     $$
Secondly, we compute $ \depth S/((I_3(P_n^2):x_{n-2}), x_{n-4})$.
$$((I_3(P_n^2):x_{n-2}), x_{n-4}) =I_3(P_{n-5}^2)+( x_{n-4}, x_{n-3}x_{n-1}, x_{n-1}x_{n}, x_{n-5}x_{n-3},x_{n-3}x_{n}).$$
$$(((I_3(P_n^2):x_{n-2}), x_{n-4}): x_{n-3})=(I_3(P_{n-5}^2)+( x_{n-4}, x_{n-3}x_{n-1}, x_{n-1}x_{n}, x_{n-5}x_{n-3},x_{n-3}x_{n})):x_{n-3}$$
$$\qquad \,\,\,\, \qquad  =I_3(P_{n-5}^2)+( x_{n-5}, x_{n-4}, x_{n-1}, x_{n})$$ $$\qquad \qquad \,  \quad=I_3(P_{n-6}^2)+( x_{n-5}, x_{n-4}, x_{n-1}, x_{n}).$$
Then by Lemma 2.2, $$\depth S/(((I_3(P_n^2):x_{n-2}), x_{n-4}): x_{n-3})= \depth S/(I_3(P_{n-6}^2)+( x_{n-5}, x_{n-4}, x_{n-1}, x_{n}))\qquad $$
$$ \qquad \qquad \qquad \qquad \qquad \qquad \qquad \qquad \quad \,\,\,  = \depth S_{n-6}/(I_3(P_{n-6}^2)+ \depth S_{{(n-6)}}/(x_{n-5}, x_{n-4}, x_{n-1}, x_{n}))$$ $$  \qquad \quad  \qquad \quad  \qquad = \depth S_{n-6}/I_3(P_{n-6}^2)+ 2.$$
While,
$$((I_3(P_n^2):x_{n-2}), x_{n-4}, x_{n-3})= I_3(P_{n-5}^2)+( x_{n-4}, x_{n-3}, x_{n-1}x_{n}).$$
By Lemma 2.2,
$$\depth S/((I_3(P_n^2):x_{n-2}), x_{n-4}, x_{n-3})= \depth S_{n-5}/I_3(P_{n-5}^2)+ \depth S_{(n-5)}/( x_{n-4}, x_{n-3}, x_{n-1}x_{n})$$ $$  \qquad \quad  \qquad \quad   = \depth S_{n-5}/I_3(P_{n-5}^2)+ 2  .$$
By induction, $$ \depth S_{n-6}/I_3(P_{n-6}^2)=\left\lceil \frac{n-6}{7} \right\rceil+ \left\lfloor \frac{n-8}{7} \right\rfloor  +1 $$ $$ \qquad \quad \qquad \quad \qquad \qquad \quad \qquad \quad \qquad \quad \,\,  \quad \leq  \left\lceil \frac{n-5}{7} \right\rceil+ \left\lfloor \frac{n-7}{7} \right\rfloor  +1 =\depth S_{n-5}/I_3(P_{n-5}^2).$$ Then this leads to $$\depth S/(((I_3(P_n^2):x_{n-2}), x_{n-4}): x_{n-3}) \leq \depth S/((I_3(P_n^2):x_{n-2}), x_{n-4}, x_{n-3}).$$ Hence, by Lemma 2.1,
$$\depth S/((I_3(P_n^2):x_{n-2}), x_{n-4}) = \depth S/(((I_3(P_n^2):x_{n-2}), x_{n-4}): x_{n-3})$$ $$ \qquad \quad \qquad \qquad \,\,\,  = \depth S_{n-6}/I_3(P_{n-6}^2)+ 2.$$
By induction, $ \depth S_{n-7}/I_3(P_{n-7}^2)\leq  \depth S_{n-6}/I_3(P_{n-6}^2)$,
This implies that $$  \depth S/((I_3(P_n^2):x_{n-2}):x_{n-4}) \leq \depth S/((I_3(P_n^2):x_{n-2}), x_{n-4}).$$ By Lemma 2.1, $$\depth S/(I_3(P_n^2):x_{n-2})= \depth S/((I_3(P_n^2):x_{n-2}):x_{n-4})=\depth S_{n-7}/I_3(P_{n-7}^2) +2.$$

Next, we compute  $\depth S/(I_3(P_n^2), x_{n-2})$. Obviously,
%$$(I_3(P_n^2), x_{n-2})=I_3(P_{n-3}^2)+(x_{n-2}, x_{n-3}x_{n-1}x_n, x_{n-4}x_{n-3}x_{n-1}, x_{n-5}x_{n-3}x_{n-1}).$$
$$(I_3(P_n^2), x_{n-2}, x_{n-3})=I_3(P_{n-4}^2)+(x_{n-2}, x_{n-3}).$$
While,
$$((I_3(P_n^2), x_{n-2}): x_{n-3})=I_3(P_{n-6}^2)+ L,\qquad \qquad $$ where, $$L=(x_{n-2}, x_{n-5}x_{n-4}, x_{n-6}x_{n-4}, x_{n-6}x_{n-5}, x_{n-7}x_{n-5}, x_{n-1}x_{n}, x_{n-4}x_{n-1}, x_{n-5}x_{n-1}).$$
Let $H=((I_3(P_n^2), x_{n-2}): x_{n-3})$. Then $$(H:x_{n-5})=I_3(P_{n-6}^2)+(x_{n-1}, x_{n-2}, x_{n-4}, x_{n-6}, x_{n-7})$$ $$ \qquad \qquad  \quad =I_3(P_{n-8}^2)+(x_{n-1}, x_{n-2}, x_{n-4}, x_{n-6}, x_{n-7}).$$
$$\qquad \qquad \quad \quad  (H, x_{n-5})=I_3(P_{n-6}^2)+(x_{n-2}, x_{n-5}, x_{n-6} x_{n-4}, x_{n-4} x_{n-1}, x_{n-1} x_{n}).$$
$$((H, x_{n-5}):x_{n-4})=I_3(P_{n-6}^2)+(x_{n-2}, x_{n-5}, x_{n-6}, x_{n-1}) \quad \quad  \quad  \quad \quad \quad $$ $$\quad \quad \quad =I_3(P_{n-7}^2)+(x_{n-2}, x_{n-5}, x_{n-6}, x_{n-1}).\,\, $$
$$((H, x_{n-5}), x_{n-4})=I_3(P_{n-6}^2)+(x_{n-2}, x_{n-4}, x_{n-5},  x_{n-1} x_{n}).\quad  \quad  \qquad  $$
By using lemma 2.1 and lemma 2.2 repeatedly, we get that $$\depth S/ (I_3(P_n^2), x_{n-2})\qquad \qquad \qquad \qquad \qquad \qquad \qquad\qquad\qquad\qquad \qquad\qquad \qquad \qquad\qquad  $$ $$\in \{\depth S/I_3(P_{n-6}^2) +2, \depth S/I_3(P_{n-7}^2)+3, \depth S/I_3(P_{n-8}^2)+3, \depth S/I_3(P_{n-4}^2)+2 \}.$$
Hence, by induction and Lemma 2.3,  $$\depth S/(I_3(P_n^2):x_{n-2})= \depth S_{n-7}/I_3(P_{n-7}^2) +2 \leq \depth S/ (I_3(P_n^2), x_{n-2}).$$
Then by Lemma 2.1 and induction, $$ \qquad \qquad \depth S/I_3(P_n^2)=\depth S/(I_3(P_n^2):x_{n-2})= \depth S_{n-7}/I_3(P_{n-7}^2) +2 $$ $$=\left\lceil \frac{n}{7} \right\rceil+ \left\lfloor \frac{n-2}{7} \right\rfloor  +1.  \qquad \,\, \,  $$ This completes the proof.

\vskip2mm

\no{\bf Corolary 2.5} \,\, Let $I_3(P_n^2)$ be  the 3-path ideal of square of a path $P_n^2$. Then the projective dimension of $S/I_3(P_n^2)$ is $$\pd({{S}}/I_3(P_{n}^{2}))=n-1-\left\lceil \frac{n}{7} \right\rceil- \left\lfloor \frac{n-2}{7} \right\rfloor.$$

\vskip2mm

\no{\bf Proof.} It follows by Auslander-Buchsbaum Formula \cite[Theorem 1.3.3]{BrH} and Theorem 2.4.

\vskip2mm

\no{\bf Proposiion 2.6} \,\, Let $I_3(P_n^2)$ be  the 3-path ideal of square of a path $P_n^2$. Then $${\text{dim}}({{S}}/I_3(P_{n}^{2}))=\left\lceil \frac{n}{4} \right\rceil+ \left\lfloor \frac{n-2}{4} \right\rfloor  +1 .$$
\vskip2mm
\no{\bf Proof.} Obvious, ${\text{dim}}({{S_3}}/I_3(P_{3}^{2}))=2$, and the result for $n=3$ holds.

Note that
$$I_3(P_n^2)=I_3(P_{n-4}^2)+ J,$$ where $$J=(x_{n-5}x_{n-4}x_{n-3}, x_{n-4}x_{n-3}x_{n-2}, x_{n-3}x_{n-2}x_{n-1}, x_{n-2}x_{n-1}x_{n}, $$
 $$\, \, \, \, \,\, \,\, \,\,\, x_{n-6}x_{n-4}x_{n-3}, x_{n-5}x_{n-3}x_{n-2}, x_{n-4}x_{n-2}x_{n-1},x_{n-3}x_{n-1}x_{n},$$
$$\, \, \, \, \,\, \,\, \,\, \, x_{n-6}x_{n-5}x_{n-3}, x_{n-5}x_{n-4}x_{n-2}, x_{n-4}x_{n-3}x_{n-1}, x_{n-3}x_{n-2}x_{n},$$
$$ \,\, \,\, \,\, \,\, \,\,\, \, x_{n-7}x_{n-5}x_{n-3}, x_{n-6}x_{n-4}x_{n-2}, x_{n-5}x_{n-3}x_{n-1}, x_{n-4}x_{n-2}x_{n}).$$

Note that $4\leq n\leq 6$. $I_3(P_{n-4}^2)=0$
%I_3(P_n^2)=I_3(P_{n-4}^2)+ J=J$. And $(x_{n-3}, x_{n-2})$ is a minimal prime ideal over $J$ and ${\text{ht}}(J)=2$. So $${\text{dim}}({{S_4}}/I_3(P_{4}^{2}))=2, {\text{dim}}({{S_5}}/I_3(P_{5}^{2}))=3, {\text{dim}}({{S_6}}/I_3(P_{6}^{2}))=4.$$
and that $(x_{n-3}, x_{n-2})$ is a minimal prime ideal over $J$ and ${\text{ht}}(J)=2$. Considering the primary decomposition of $I_3(P_{n}^{2})$. For any minimal prime $p$ over $I_3(P_{n-4}^{2})$, $p+(x_{n-3}, x_{n-2})$ is a  minimal prime over $I_3(P_{n}^{2})$. Hence, $${\text{ht}} (I_3(P_{n}^{2}))\leq {\text{ht}}(I_3(P_{n-4}^{2}))+2.$$ On the other hand, any one of $\{ x_{n-3}, x_{n-2}, x_{n-1}, x_{n}  \}$ can't cover all generaters of $J$, so $${\text{ht}}(I_3(P_{n}^{2}))\geq {\text{ht}}(I_3(P_{n-4}^{2}))+2.$$ Then $${\text{ht}}(I_3(P_{n}^{2}))= {\text{ht}}(I_3(P_{n-4}^{2}))+2,$$  and by \cite[Corollary 3.1.7]{Vi}, $${\text{dim}}({{S}}/I_3(P_{n}^{2}))= {\text{dim}}{{S}}- {\text{ht}}(I_3(P_{n}^{2}))\qquad \qquad  $$ $$\qquad \qquad \qquad \quad= {\text{dim}}{{S}_{n-4}}- {\text{ht}}(I_3(P_{n-4}^{2}))+2$$ $$\qquad \qquad \qquad ={\text{dim}}({{S}_{n-4}}/I_3(P_{n-4}^{2}))+2.$$
By induction, $${\text{dim}}({{S}_{n}}/I_3(P_{n}^{2}))=\left\lceil \frac{n}{4} \right\rceil+ \left\lfloor \frac{n-2}{4} \right\rfloor  +1.$$ \vskip2mm

By Theorem  2.4 and Proposition 2.6, if ${\text{dim}}({{S}}/I_3(P_{n}^{2}))= \depth ({{S}_{n}}/I_3(P_{n}^{2}))$, then this leads to $n=3, 4.$ In fact, $I_3(P_{3}^{2})=(x_1x_2x_3)$, this is a complete intersection, and $S/I_3(P_{3}^{2})$ is Cohen-Macaulay. While for $I_3(P_{4}^{2})=(x_1x_2x_3, x_2x_3x_4, x_1x_3x_4, x_1x_2x_4)$, we regard $I_3(P_{4}^{2})$ as the stanley-Reiner ideal of a 1-skeleton $\gamma$ of a simplex $\delta$ on $\{1,2,3,4\}$. This simplex $\delta$ is shellabe and then by \cite[Theorem 8.2.18]{HH},  the 1-skeleton $\gamma$ is pure and shellable. So   $S/I_3(P_{4}^{2})$ is Cohen-Macaulay by \cite[Theorem 8.2.6]{HH}.

\vskip2mm

\no{\bf Theorem 2.7}  Let $I_3(P_n^2)$ be  the 3-path ideal of square of a path $P_n^2$.  Then $S/I_3(P_{n}^{2})$ is a Cohen-Macaulay ideal if and only if $n=3, 4.$

\vskip2mm

Let $R$ be a polynomial ring over a field $K$ and $I$ a homogeneous ideal in $R$. Brodmann \cite{Br} showed that $\depth R/I^n$ is a constant for sufficiently large $n$.
In the following, we consider the limit behavior of depth of powers of 3-path ideal of square of a path graph. \vskip2mm

\no{\bf Proposition 2.8}\,  Let $I_3(P_n^2)$ be  the 3-path ideal of square of a path $P_n^2$. Then

(1)  $\depth({{S}_{3}}/(I_3(P_{3}^{2}))^t)=2$ for any $t>0 $.

(2) $\depth({{S}_{4}}/(I_3(P_{4}^{2}))^t)=0$ for any $t> 2$.

(3) $\depth({{S}_{n}}/(I_3(P_{n}^{2}))^t)=0$ for any $n\geq 5$ and $t>\left\lfloor \frac{n-2}{3} \right\rfloor$. \vskip2mm

\no {\bf Proof.} Obviously, $I_3(P_{3}^{2})=(x_1x_2x_3)$ and $\depth({{S}_{3}}/(I_3(P_{3}^{2}))^t)=2$ for any $t>0 $.

By Theorem 2.6,  $I_3(P_{4}^{2})$ is Cohen-Macaulay. In particulay, $${\text{dim}}({{S}_{4}}/(I_3(P_{4}^{2})))= \depth({{S}_{4}}/(I_3(P_{4}^{2})))=2.$$
If ${{S}_{4}}/(I_3(P_{4}^{2}))^2$ is Cohen-Maucaulay, then $${\text{Ass}}({{S}_{4}}/(I_3(P_{4}^{2}))^2)\subseteq {\text{Ass}}({{S}_{4}}/(I_3(P_{4}^{2}))).$$ Therefore,  $$(I_3(P_{4}^{2}))^2=(I_3(P_{4}^{2}))^{(2)}.$$  Note that $$I_3(P_{4}^{2})=(x_1, x_2)\cap (x_1, x_3)\cap (x_1, x_4)\cap(x_2, x_3)\cap (x_2, x_4)\cap (x_3, x_4).$$   Obviously, $x_1^2x_2x_3x_4\in (I_3(P_{4}^{2}))^{(2)}\setminus (I_3(P_{4}^{2}))^{2}$. This is a contradiction.
%Note that $(I_3(P_{4}^{2})=(x_1, x_2)\cap (x_1, x_3)\cap (x_1, x_4)\cap(x_2, x_3)\cap (x_2, x_4)\cap (x_3, x_4)$.
Hence, $$\depth({{S}_{4}}/(I_3(P_{4}^{2}))^2)\leq 1.$$
Suppose that $\mathfrak{m}=(x_1, x_2, x_3, x_4)\in {\text{Ass}}{{S}_{4}}/(I_3(P_{4}^{2}))^2$. Then there exists $$u\in ((I_3(P_{4}^{2})^2: \mathfrak{m})\setminus (I_3(P_{4}^{2}))^2 ,\, \,  {\text{deg}}(u)\geq 5.$$ %Let $u=x_{i_1}x_{i_2}x_{i_3}x_{i_4}x_{i_5}x_{i_6}v, i_1\leq i_2\leq i_3\leq i_4\leq i_5 \leq i_6 \, v $ a monomial.
Since $u\notin (I_3(P_{4}^{2}))^2$, $x_ju\in (I_3(P_{4}^{2}))^2$ for any $j$ and $I_3(P_{4}^{2})$ is squarefree,  $\text{deg}_{x_j}u< 2$ for any $j$, that's,  the degree of each $x_j$ appearing  in $u$ is less than 2. This is impossible, since  there are only four variables appearing in $u$.
This shows that $\depth({{S}_{4}}/(I_3(P_{4}^{2}))^2)\geq 1$. Hence, $$\depth({{S}_{4}}/(I_3(P_{4}^{2})))=1.$$
Note that $w=x_1^2x_2^2x_3^2x_4^2\in (I_3(P_{4}^{2})^3: \mathfrak{m})\setminus I_3(P_{4}^{2})^3 $. Then,
 $$wv\in (I_3(P_{4}^{2}))^{3+k}: \mathfrak{m})\setminus (I_3(P_{4}^{2}))^{3+k}$$ for any generator $v\in G((I_3(P_{4}^{2}))^k)$. This shows that $\mathfrak{m}\in {\text{Ass}} S_4/(I_3(P_{4}^{2}))^{t} $ and $$\depth({{S}_{4}}/(I_3(P_{n}^{2}))^t)=0$$ for any $t> 2$.

Assume that $n\geq 5$. Let $t_0=\left\lfloor \frac{n-2}{3} \right\rfloor +1$.

Take $a=x_{r+1}x_{r+2}\cdots x_{n}$, where $r\equiv n-2 \pmod{3}, \, 0\leq r\leq 2.$ Note that
$$3t_0-(n-r)=3(\left\lfloor \frac{n-2}{3} \right\rfloor +1)-(n-r)=n-2-r+3-(n-r)=1.$$ Then $ {\text{deg}}(a)=n-r=3t_0-1$ and $a\notin I^{t_0}$.
In the following, we will show that $$a\in (I^{t_0}: \mathfrak{m}),\, \text{that's,\, for\, any} \, \, l, \, ax_l\in I^{t_0}.$$
In case that $n=3k+2$: \, $r=0$, $a=x_{1}x_{2}\cdots x_{n}$.

If $l\equiv 0\pmod{3}$. Let $l=3s$, $$ax_{3s}=(x_1x_2x_3)\cdots(x_{3s-2}x_{3s-1}x_{3s})(x_{3s}x_{3s+1}x_{3s+2})\cdots (x_{3k}x_{3k+1}x_{3k+2})\in I^t.$$

If $l\equiv 1\pmod{3}$. Let $l=3s+1$, $$ \quad \quad \,  ax_{3s+1}=\cdots (x_{3s+1}x_{3s+2}x_{3(s+1)+1})(x_{3s+1}x_{3(s+1)}x_{3(s+1)+2})\cdots (x_{3k}x_{3k+1}x_{3k+2}) \in I^t.$$

If $l\equiv 2\pmod{3}$. Let $l=3s+2$, $$\quad \quad \,  ax_{3s+2}=\cdots (x_{3s+1}x_{3s+2}x_{3(s+1)})(x_{3s+2}x_{3(s+1)+1}x_{3(s+1)+2})\cdots (x_{3k}x_{3k+1}x_{3k+2})\in I^t.$$
In case that $n=3k+3$:\, $r=1$, $a=x_{2}x_{3}\cdots x_{n}$.

If $l\equiv 0\pmod{3}$. Let $l=3s$, $$  \quad ax_{3s}=(x_2x_3x_4)\cdots(x_{3s-1}x_{3s}x_{3s+1})(x_{3s}x_{3s+2}x_{3(s+1)})\cdots (x_{3k+1}x_{3k+2}x_{3k+3})\in I^t.$$

If $l\equiv 1\pmod{3}$. Let $l=3s+1$, $$    ax_{3s+1}=\cdots (x_{3s-1}x_{3s}x_{3s+1})(x_{3s+1}x_{3s+2}x_{3(s+1)})\cdots (x_{3k+1}x_{3k+2}x_{3k+3})\in I^t.$$

If $l\equiv 2\pmod{3}$. Let $l=3s+2$, $$  ax_{3s+2}=\cdots (x_{3s-1}x_{3s}x_{3s+2})(x_{3s+1}x_{3s+2}x_{3(s+1)})\cdots (x_{3k+1}x_{3k+2}x_{3k+3})\in I^t.$$
In case that $n=3k+4$: \, $r=2$, $a=x_{3}x_{4}\cdots x_{n}$.

If $l\equiv 0\pmod{3}$. Let $l=3s$,$$ \quad \quad ax_{3s}=(x_3x_4x_5)\cdots(x_{3(s-1)}x_{3s-2}x_{3s})(x_{3s-1}x_{3s}x_{3s+1})\cdots (x_{3k+2}x_{3k+3}x_{3k+4})\in I^t.$$

If $l\equiv 1\pmod{3}$. Let $l=3s+1$, $$\quad  ax_{3s+1}=\cdots (x_{3(s-1)}x_{3s-1}x_{3s+1})(x_{3s-2}x_{3s}x_{3s+1})\cdots (x_{3k+2}x_{3k+3}x_{3k+4})\in I^t.$$

If $l\equiv 2\pmod{3}$. Let $l=3s+2$, $$\quad \quad \, ax_{3s+2}=\cdots (x_{3s}x_{3s+1}x_{3s+2})(x_{3s+2}x_{3(s+1)}x_{3(s+1)+1})\cdots (x_{3k+2}x_{3k+3}x_{3k+4})\in I^t.$$
Now we have  shown that $a\in (I^{t_0}: \mathfrak{m})\setminus I^{t_0}$. Furtherly, it is clear that  $au\in (I^{t_0+j}: \mathfrak{m})\setminus I^{t_0+j}$ for any generator $u\in G(I_3(P_{n}^{2})^j)$. Therefore, $\mathfrak{m}\in {\text{Ass}} S/I^{t}$ and $\depth({{S}_{3}}/(I_3(P_{n}^{2}))^t)=0$ for any $t\geq t_0$.

\vskip2mm

\noindent
{\bf Acknowledgments.} We gratefully acknowledge the computer algebra
system {\tt CoCoA} \cite{Co} which was invaluable in our work on this paper.

%\bibitem{LT} Le Xuan Dung, Truong Thi Hien, Hop D. Nguyen£¬ Tran
%Nam Trung, Regularity and Koszul property of symbilic powers of
%monomial ideals, Math. Z. 298(2021), no3-4,1487-1522.

%John A. Eagon, Victor Reiner, Resolutions of Stanley-Reisner rings and Alexander duality. J. Pure
%Appl. Algebra 130(1998), no. 3, 265¨C275.

% J¡§urgen Herzog, Takayuki Hibi, Componentwise linear ideals. Nagoya Math. J. 153 (1999), 141¨C153.

\vskip5mm

Liuqing Yang

Department of Mathematics, Soochow University, Suzhou,
215006,China

E-mail: 20214007001@stu.suda.edu.cn

\vskip3mm

Lizhong Chu

Department of Mathematics, Soochow University, Suzhou,
215006,China

E-mail: chulizhong@suda.edu.cn

\end{document}